\documentclass[12pt]{article}
\usepackage[pctex32]{graphics}
\usepackage{amsthm,amsmath,amssymb,amscd,verbatim,epsfig}
\usepackage{amssymb}
\usepackage{graphicx}
\usepackage{epsfig}
\usepackage{graphics,pifont}
\newcommand{\beq}{\begin{equation}}
\newcommand{\eeq}{\end{equation}}

\date{}

\newcommand{\ld}{\lambda}

\newcommand{\f}{\frac}
\newcommand{\ra}{\rightarrow}

\begin{document}

\title{The\ Minimal\ Number\ of\ Periodic\ Orbits\ of\ Periods\ Guaranteed\ in\ Sharkovskii's\ Theorem}
\author{Bau-Sen Du \\ [.5cm]
Institute of Mathematics \\
Academia Sinica \\
Taipei 11529, Taiwan \\
dubs@math.sinica.edu.tw \\ [.2cm]
(Bull. Austral. Math. Soc. 31(1985), 89-103. Corrigendum: 32 (1985), 159) \\}
\maketitle
\begin{abstract}
Let $f(x)$ be a continuous function from a compact real interval into itself with a periodic orbit of minimal period $m$, where $m$ is not an integral power of 2.  Then, by Sharkovskii's theorem, for every positive integer $n$ with $m \ra n$ in the Sharkovskii's ordering defined below, a lower bound on the number of periodic orbits of $f(x)$ with minimal period $n$ is 1.  Could we improve this lower bound from 1 to some larger number?  In this paper, we give a complete answer to this question.
\end{abstract}

%37E05, (37C25, 37E15)

\section{Introduction}

Let $I$ be a compact real interval and let $f \in C^0(I, I)$.  For any $x_0$ in $I$ and any positive integer $k$, we let $f^k(x_0)$ denote the $k$th iterate of $x_0$ under $f$ and call $\{f^k(x_0) : k \ge 0 \}$ the orbit of $x_0$ (under $f$).  If $f^m(x_0) = x_0$ for some positive integer $m$, we call $x_0$ a periodic point of $f$ and call the cardinality of the orbit of $x_0$ (under $f$) the minimal period of $x_0$ and of the orbit (under $f$).  If $f$ has a periodic orbit of a period $m$, must $f$ also have periodic orbits of periods $n \ne m$ ?  In 1964, Sharkovskii {\bf{\cite{sh}}} (see {\bf{\cite{bl, bu, ho, li, st, str}}} also)had given a complete answer to this question.  Arrange the positive integers according as the following new order (called Sharkovskii ordering):
$$
3 \prec 5 \prec 7 \prec  9 \prec \cdots \prec 2 \cdot 3 \prec 2 \cdot 5 \prec 2 \cdot 7 \prec  2 \cdot 9 \prec \cdots \prec 2^2 \cdot 3 \prec 2^2 \cdot 5 \prec 2^2 \cdot 7 \prec  2^2 \cdot 9 \prec \cdots$$ $$\prec \cdots \prec 2^3 \prec 2^2 \prec 2 \prec 1.
$$

\noindent
Sharkovskii's theorem says that any function $f \in C^0(I, I)$ with a periodic orbit of minimal period $m$ must also have at least one periodic orbit of minimal period $n$ precisely when $m \prec n$ in the above Sharkovskii ordering.  Therefore, for every positive integer $n$ with $m \prec n$, the number 1 is a lower bound on the number of distinct periodic orbits of $f$ with minimal periodi $n$.  One question arises naturally: Could we improve this lower bound from 1 to some larger number ?

In 1976 Bowen and Franks {\bf{\cite{bo}}} showed, among other things, that if $f \in C^0(I, I)$ has a periodic orbit of minimal period $n = 2^dm$, where $m > 1$ is odd, then there is a number $M_n$ (independent of $f$) such that, for all integers $k \ge M_n$, $f$ has at least $(2^{k/m})/(2^dk)$ distinct periodic orbits of minimal period $2^dk$.  

In 1979, Jonker {\bf{\cite{jo}}} also obtained a similar result on a class of unimodal maps.  If $c$ is an interior point of $I$, let $S_c$ denote the collection of all $f \in C^0(I, I)$ which has either one maximum or one minimum point at $c$, and is strictly monotone on each component of $I - \{ c \}$ with $f(\partial I) \subset \partial I$.  Jonker showed, among other things, that if $m$, $n$ are any two odd integers with $1 < m < n$, and if $f \in S_c$ has a periodic orbit of minimal period $2^km$, where $k \ge 0$ is any integer, then $f$ must also have at least $2^{(n-m)/2}$ distinct periodic orbits of minimal period $2^kn$.  

In {\bf{\cite{du3}}}, a result along this line is also obtained.  However, that result is only a partial one.  In this paper we give a complete answer to that question.  

In Section 2 we state our main results (Theorems 1, 2, and 3).  In Section 3 we describe the method used to prove them.  This method is the same as that used in {\bf{\cite{du2,du3}}}.  The proofs of Theorems 1 and 2 will appear in Sections 4 and 5.  Theorem 3 then follows easily from Theorems 1 and 2.

\section{Statement of main results}

Let $\phi(m)$ be an integer-valued function defined on the set of all positive integers.  If $m=p_1^{k_1}p_2^{k_2} \cdots p_r^{k_r}$, where the $p_i$'s are distinct prime numbers, $r$ and $k_i$'s are positive integers, we define $$\Phi(1, \phi)=\phi(1)$$ and $\Phi(m, \phi) =$
$$
\phi(m)-\sum_{i=1}^r \phi(\f m{p_i})+\sum_{i_1<i_2} \phi(\f m{p_{i_1}p_{i_2}})
- \sum_{i_1<i_2<i_3} \phi(\f m{p_{i_1}p_{i_2}p_{i_3}}) + \cdots 
+ (-1)^r \phi(\f m{p_1p_2 \cdots p_r}),
$$
\noindent
where the summation $\sum_{i_1<i_2< \cdots < i_j}$ is taken over all integers $i_1, i_2, \cdots, i_j$ with $1 \le i_1 < i_2 < \cdots < i_j \le r$.  If, when considered as a sequence, $<\phi(m)>$ is the Lucas sequence, that is, is $\phi(1) = 1, \phi(2) = 3$, and $\phi(m+2) = \phi(m+1) + \phi(m)$ for all positive integers $m$, then, for simplicity, we denote $\Phi(m, \phi)$ as $\Phi_1(m)$.  Note that, if $f \in C^0(I, I)$ and if, for every positive integer $m$, $\phi(m)$ is the number of distinct solutions of the equation $f^m(x) = x$, then $\Phi(m, \phi)$ is, by the standard inclusion-exclusion argument, the number of periodic points of $f$ with minimal period $m$.  Now we can state the following theorem.  

\noindent
{\bf Theorem 1.}
{\it Let $f : [1, 3] \ra [1, 3]$ be defined by $f(x) = -2x + 5$ if $1 \le x \le 2$ and $f(x) = x - 1$ if $2 \le x \le 3$.  Then the following hold:
\begin{itemize}

\item[(a)] 
for every positive integer $m$, if $a_m$ is the number of distinct solutions of the equation $f^m(x) = x$, then the sequence $<a_m>$ is the Lucas sequence;

\item[(b)]
for every positive integer $m$, $f$ has exactly $\Phi_1(m)/m$ distinct periodic orbits of minimal period $m$;

\item[(c)]
the sequence $<\Phi_1(m)/m>$ is strictly increasing for $m \ge 6$ and $\lim_{n \to \infty} [\Phi_1(m+1)/(m+1)]/[\Phi_1(m)/m] = (1 + \sqrt{5})/2$.
\end{itemize}}

Fix any integer $n > 1$ and let $$Q_n = \{(1, n+1) \} \cup \{ (m, 2n+2-m) : 2 \le m \le n \} \cup \{ (m, 2n+1-m) : n+1 \le m \le 2n \}.$$
\noindent
For all integers $i, j$, and $k$, with $1 \le i, j \le 2n$ and $ k \ge 1$, we define $b_{k,i,j,n}$ recursively as follows:
$$
b_{1,i,j,n} =
\begin{cases}
1, & \text{if} \,\,\, (i,j) \in Q_n, \\
0, & \text{otherwise}, \\
\end{cases}
$$

\noindent
and

$$
b_{k+1,i,j,n} = 
\begin{cases}
b_{k,i,2n+1-j,n} + b_{k,i,n+1,n}, & \text{if} \,\,\, 1 \le j \le n-1, \\
b_{k,i,n,n} + b_{k,i,n+1,n}, & \text{if} \,\,\, j = n, \\
b_{k,i,1,n}, & \text{if} \,\,\, j = n+1, \\
b_{k,i,2n+2-j,n}, & \text{if} \,\,\, n+2 \le j \le 2n. \\
\end{cases}
$$

\noindent
We also define $c_{k,n}$ by letting $$c_{k,n} = \sum_{i=1}^{2n} b_{k,i,i,n} + b_{k,n+1,n,n} + \sum_{i=n+2}^{2n} b_{k,i,n+1,n}.$$  Note that these sequences $< b_{k,i,j,n} >$ and $< c_{k,n} >$ have the following six properties.  Some of these will be used later in the proofs of our main results.  (Recall that $n > 1$ is fixed.)  

(i) The sequence $<b_{k,1,n,n}>$ is increasing, and for all integers $k \ge 2$, we have $b_{k,1,n,n} \ge b_{k,n+1,n,n}$ and $b_{k,1,i+1,n} \ge b_{k,1,i,n}$ for all $1 \le i \le n-1$.

(ii) The sequences $<b_{k,1,j,n}>$, $1 \le j \le n$, and $<b_{k,n+1,n,n}>$ can also be obtained by the following recursive formulas:
$$
b_{1,1,j,n} = 0, \,\,\, 1 \le j \le n,
$$
$$
b_{2,1,j,n} = 1, \,\,\, 1 \le j \le n,
$$
$$
b_{1,n+1,n,n} = b_{2,n+1,n,n} = 1,
$$
$$
b_{1,n+1,j,n} = b_{2,n+1,j,n} = 0, \,\,\, 1 \le j \le n-1.
$$

For $i = 1$ or $n+1$, and $k \ge 1$,
$$
b_{k+2,i,n,n} = b_{k,i,1,n} + b_{k+1,i,n,n},
$$
$$
b_{k+2,i,j,n} = b_{k,i,1,n} + b_{k,i,j+1,n}, \,\,\, 1 \le j \le n-1.
$$

(iii) For every positive integer $k$, $c_{k+2n-2,n}$ can also be obtained by the following formulas:
$$
c_{k+2n-2,n} = b_{k+2n-2,n+1,n,n} + 2 \sum_{j=1}^n b_{k+2n-2j,1,j,n} \qquad\qquad\qquad\qquad$$ $$\quad\qquad\, = b_{k+2n-2,n+1,n,n} + 2nb_{k,1,n,n} + \sum_{i=2}^n (2^i-2)b_{k,1,n+1-i,n}.
$$

\noindent
The first identity also holds for all integers $k$ with $-2n+3 \le k \le 0$ provided we define $b_{k,1,j,n} = 0$ for all $-2n+3 \le k \le 0$ and $1 \le j \le n$.  

(iv) For all integers $k$ with $1 \le k \le 2n$, $c_{2k,n} = 2^{k+1} - 1$.

(v) For all integers $k$ with $n+1 \le k \le 3n$, $c_{2k+1,n} = 2c_{2k+1,n+1} - 1$.

(vi) Since, for every positive integer $k \ge 2n+1$, $$b_{k,1,n,n} = b_{k-1,1,n,n} + \sum_{i=2}^{2n} (-1)^ib_{k-i,1,n,n},$$ there exist $2n+1$ nonzero constants $\alpha_j$'s such that $b_{k,1,n,n} = \sum_{j=1}^{2n+1} \alpha_jx_j^k$ for all positive integers $k$, where $\{ x_j : 1 \le j \le 2n+ 1 \}$ is the set of all zeros (including complex zeros) of the polynomial $x^{2n+1} - 2x^{2n-1} - 1$.  

For all positive integers $k$, $m$, $n$, with $n > 1$, we let $\phi_n(k) = c_{k,n}$ and let $\Phi_n(m) = \Phi(m, \phi_n)$, where $\Phi$ is defined as above.  Now we can state the following theorem.

\noindent
{\bf Theorem 2.}
{\it For every integer $n > 1$, let $f_n : [1, 2n+1] \ra [1, 2n+1]$ be the continuous function with the following six properties:
\begin{itemize}

\item[(1)]
$f_n(1) = n + 1$,

\item[(2)]
$f_n(2) = 2n + 1$,

\item[(3)]
$f_n(n+1) = n + 2$,

\item[(4)]
$f_n(n+2) = n$, 

\item[(5)]
$f_n(2n+1) = 1$, and 

\item[(6)]
$f_n$ is linear on each component of the complement of the set $\{ 2, n+1, n+2 \}$ in $[1, 2n+1]$.
\end{itemize}

\noindent
Then the following hold:
\begin{itemize}

\item[(a)]
For every positive integer $k$, the equation $f_n^k(x) = x$ has exactly $c_{k,n}$ distinct solutions;

\item[(b)]
For every positive integer $m$, $f_n(x)$ has exactly $\Phi_n(m)/m$ distinct periodic orbits of minimal period $m$;

\item[(c)]
$\lim_{m \to \infty} (\log [\Phi_n(m)/m])/m = \log \ld_n$, where $\ld_n$ is the (unique) positive (and the largest in absolute value) zero of the polynomial $x^{2n+1} - 2x^{2n-1} - 1$.  
\end{itemize}}

From Theorems 1 and 2 above and Theorem 2 of {\bf [12}, p. 243 {\bf ]}, we easily obtain the following result.  

\noindent
{\bf Theorem 3.}
{\it Assume that $f \in C^0(I, I)$ has a periodic orbit of minimal period $s = 2^k(2n+1)$, where $n \ge 1$ and $k \ge 0$, and no periodic orbits of minimal period $r$ with $r \prec s$ in the Sharkovskii ordering.  Then for every positive integer $t$ with $s \prec t$ in the Sharkovskii ordering, $f$ has at least $\Phi_n(t/2^k)/(t/2^k)$ (sharp) distinct periodic orbits of minimal period $t$.}

\noindent
{\bf Remark 1.}
We call attention to the fact that there exist continuous functions from $I$ into $I$ with exactly one periodic orbit of minimal period $2^i$ for every positive integer $i$ (and two fixed points), but no other periods (see {\bf{\cite{mi}}}).

\noindent
{\bf Remark 2.} 
With the help of Theorem 2 of {\bf [12}, p. 243 {\bf ]} on the distribution along the real line of points in a periodic orbit of odd period $n > 1$, when there are no periodic orbits of odd period $m$ with $1 < m < n$, our results give a new proof of Sharkovskii's theorem.  

\noindent
{\bf Remark 3.}
Table 1 (see next page) lists the first 31 values of $\Phi_n(m)/m$ for $1 \le n \le 5$.  It seems that, for all positive integers $n$ and $m$, we have $$\Phi_n(2m+1)/(2m+1) = 2^{m-n} \,\,\, \text{for} \,\,\, n \le m \le 3n+1,$$ and $$\Phi_n(2m+1)/(2m+1) > 2^{m-n} \,\,\, \text{for} \,\,\, n > 3n+1.$$

\noindent
{\bf Remark 4.}
For all positive integers $k$ and $m$, let $\psi(k) = 2^k$ and $\Psi(m) = \Phi(m, \psi)$, where $\Phi$ is defined as in Section 2.  It is obvious that $\Psi(m)/m$ is the number of distinct periodic orbits of minimal period $m$ for, say, the mapping $g(x) = 4x(1-x)$ from $[0, 1]$ onto itself.  Since, for all positive integers $k$ and $n$ with $1 \le k \le 2n$, $c_{2k,n} = 2^{k+1} - 1$ and $c_{1,n} = 1$, we obtain that $\Phi_n(2k+2)/(2k+2) = \Psi(k+1)/(k+1)$ for all $1 \le k \le 2n$.    It seems that $\Phi_n(2k+2)/(2k+2) > \Psi(k+1)/(k+1)$ for all $k > 2n$.  But note that $$\lim_{k \to \infty} (\log [\Phi_n(2k+2)/(2k+2)])/(2k+2) = \log \ld_n > \f 12 \log 2 = \f 12 \lim_{k \to \infty} (\log [\Psi(k+1)/(k+1)])/(k+1),$$ where $\ld_n$ is the unique positive zero of the polynomial $x^{2n+1} - 2x^{2n-1} - 1$.

$$\text{Table 1}$$
\medskip
\noindent
$m$ \quad\qquad $\Phi_1(m)/m$ \qquad $\Phi_2(m)/m$ \qquad $\Phi_3(m)/m$ \qquad $\Phi_4(m)/m$ \qquad $\Phi_5(m)/m$ \qquad $\Psi(m)/m$
1 \quad \qquad\qquad 1 \qquad\qquad\qquad 1 \qquad\qquad\qquad 1 \quad\qquad\qquad 1 \quad\quad\qquad\qquad 1 \quad\qquad\qquad 2 \\
2 \quad\qquad\qquad 1 \qquad\qquad\qquad 1 \qquad\qquad\qquad 1 \quad\qquad\qquad 1 \quad\quad\qquad\qquad 1 \quad\qquad\qquad 1 \\  
3 \quad\qquad\qquad 1 \qquad\qquad\qquad 0 \qquad\qquad\qquad 0 \quad\qquad\qquad 0 \quad\quad\qquad\qquad 0 \quad\qquad\qquad 2 \\
4 \quad\qquad\qquad 1 \qquad\qquad\qquad 1 \qquad\qquad\qquad 1 \quad\qquad\qquad 1 \quad\quad\qquad\qquad 1 \quad\qquad\qquad 3 \\
5 \quad\qquad\qquad 2 \qquad\qquad\qquad 1 \qquad\qquad\qquad 0 \quad\qquad\qquad 0 \quad\quad\qquad\qquad 0 \quad\qquad\qquad 6 \\
6 \quad\qquad\qquad 2 \qquad\qquad\qquad 2 \qquad\qquad\qquad 2 \quad\qquad\qquad 2 \quad\quad\qquad\qquad 2 \quad\qquad\qquad 9 \\
7 \quad\qquad\qquad 4 \qquad\qquad\qquad 2 \qquad\qquad\qquad 1 \quad\qquad\qquad 0 \quad\quad\qquad\qquad 0 \quad\qquad\qquad 18 \\
8 \quad\qquad\qquad 5 \qquad\qquad\qquad 3 \qquad\qquad\qquad 3 \quad\qquad\qquad 3 \quad\quad\qquad\qquad 3 \quad\qquad\qquad 30 \\
9 \quad\qquad\qquad 8 \qquad\qquad\qquad 4 \qquad\qquad\qquad 2 \quad\qquad\qquad 1 \quad\quad\qquad\qquad 0 \quad\qquad\qquad 56 \\
10 \qquad\quad\quad 11 \qquad\qquad\qquad 6 \qquad\qquad\qquad 6 \quad\qquad\qquad 6 \quad\quad\qquad\qquad 6 \quad\qquad\qquad 99 \\
11 \qquad\quad\quad 18 \qquad\qquad\qquad 8 \qquad\qquad\qquad 4 \quad\qquad\qquad 2 \quad\quad\qquad\qquad 1 \quad\qquad\qquad 186 \\
12 \qquad\quad\quad 25 \qquad\qquad\qquad 11 \qquad\qquad\qquad 9 \quad\qquad\qquad 9 \quad\quad\qquad\qquad 9 \,\,\qquad\qquad 335 \\
13 \qquad\quad\quad 40 \qquad\qquad\qquad 16 \qquad\qquad\qquad 8 \quad\qquad\qquad 4 \quad\quad\qquad\qquad 2 \,\,\qquad\qquad 630 \\
14 \qquad\quad\quad 58 \qquad\qquad\qquad 23 \qquad\qquad\qquad 18 \qquad\qquad 18 \quad\quad\qquad\qquad 18 \qquad\qquad 1161 \\
15 \qquad\quad\quad 90 \qquad\qquad\qquad 32 \qquad\qquad\qquad 16 \quad\qquad\qquad 8 \quad\quad\qquad\qquad 4 \qquad\qquad 2182 \\
16 \qquad\quad\quad 135 \qquad\qquad\qquad 46 \qquad\qquad\qquad 32 \quad\qquad\qquad 30 \quad\qquad\qquad 30 \qquad\qquad 4080 \\
17 \qquad\quad\quad 210 \qquad\qquad\qquad 66 \quad\qquad\qquad\,\, 32 \quad\qquad\qquad 16 \quad\qquad\qquad 8 \qquad\qquad 7710 \\
18 \qquad\quad\quad 316 \qquad\qquad\qquad 94 \qquad\qquad\qquad 61 \quad\qquad\qquad 56 \quad\qquad\qquad 56 \qquad\qquad 14560 \\
19 \qquad\quad\quad 492 \qquad\qquad\qquad 136 \quad\qquad\qquad 64 \quad\qquad\qquad 32 \quad\qquad\qquad 16 \qquad\qquad 27594 \\
20 \qquad\quad\quad 750 \qquad\qquad\qquad 195 \qquad\qquad\,\,\, 115 \quad\qquad\qquad 101 \quad\qquad\qquad 99 \qquad\qquad 52377 \\
21 \qquad\quad\quad 1164 \quad\qquad\qquad\,\, 282 \quad\qquad\qquad 128 \quad\qquad\qquad 64 \quad\qquad\qquad 32 \qquad\qquad 99858 \\
22 \qquad\quad\quad 1791 \quad\qquad\qquad\,\, 408 \quad\qquad\qquad 224 \qquad\qquad 191 \quad\qquad\qquad 186 \qquad\qquad 190557 \\
23 \qquad\quad\quad 2786 \quad\qquad\qquad\,\, 592 \quad\qquad\qquad 258 \qquad\qquad 128 \quad\qquad\qquad\,\, 64 \qquad\qquad 364722 \\
24 \qquad\quad\quad 4305 \quad\qquad\qquad\,\, 856 \quad\qquad\qquad 431 \qquad\qquad 351 \quad\qquad\qquad 337 \qquad\qquad 698870 \\
25 \qquad\quad\quad\,\, 6710 \quad\qquad\qquad 1248 \qquad\qquad\,\, 520 \qquad\qquad 256 \quad\qquad\qquad 128 \qquad\qquad 1342176 \\
26 \qquad\quad\quad 10420 \,\,\,\qquad\qquad\,\, 1814 \qquad\qquad\,\, 850 \qquad\qquad\,\, 668 \qquad\qquad\,\, 635 \qquad\qquad 2580795 \\
27 \qquad\quad\quad 16264 \quad\qquad\qquad 2646 \qquad\qquad 1050 \qquad\qquad\,\, 512 \qquad\qquad\,\, 256 \qquad\qquad 4971008 \\
28 \qquad\quad\quad 25350 \quad\qquad\qquad 3858 \qquad\qquad 1673 \qquad\qquad 1257 \qquad\qquad 1177 \qquad\qquad\,\, 9586395 \\
29 \qquad\quad\quad 39650 \quad\qquad\qquad 5644 \qquad\qquad 2128 \qquad\qquad 1026 \qquad\qquad\,\, 512 \qquad\qquad 18512790 \\
30 \qquad\quad\quad 61967 \quad\qquad\qquad 8246 \qquad\qquad 3328 \qquad\qquad 2402 \qquad\qquad 2220 \qquad\qquad 35790267 \\
31 \qquad\quad\quad 97108 \qquad\qquad\,\, 12088 \qquad\qquad 4320 \qquad\qquad 2056 \qquad\qquad 1024 \qquad\qquad 69273666 \\

\section{Symbolic representation for continuous piecewise linear functions}

In this section we describe a method.  This method was first introduced in {\bf{\cite{du1}}}, and then generalized in {\bf{\cite{du2}}} to construct, for every positive integer $n$, a continuous piecewise linear function from $[0, 1]$ into itself which has a periodic orbit of minimal period 3, but with the property that almost all (in the sense of Lebesgue) points of $[0, 1]$ are eventually periodic of minimal period $n$ with the periodic orbit the same as the orbit of a fixed known period $n$ point.  The same method was also used in {\bf{\cite{du3}}} to give a new proof of a result of Block {\it et al} {\bf{\cite{bl}}} on the topological entropy of interval maps.  In this paper we will use this method to prove our main results.  

Throughout this section, let $g$ be a continuous piecewise linear function from the interval $[c, d]$ into itself.  We call the set $\{(x_i, y_i) : i = 1,2, \cdots, k \}$ a set of nodes for (the graph of) $y = g(x)$ if the following three conditions hold:
\begin{itemize}

\item[(1)]
$k \ge 2$,

\item[(2)]
$x_1 = c$, $x_k = d$, $x_1 < x_2 < \cdots < x_k$, and 

\item[(3)]
$g$ is linear on $[x_i, x_{i+1}]$ for all $1 \le i \le k-1$ and $y_i = g(x_i)$ for all $1 \le i \le k$.
\end{itemize}

\noindent
For any such set, we will use its $y$-coordinates $y_1, y_2, \cdots, y_k$ to represent the graph of $y = g(x)$ and call $y_1y_2 \cdots y_k$ (in that order) a (symbolic) representation for (the graph) of $y = g(x)$.  For $1 \le i < j \le k$, we will call $y_iy_{i+1} \cdots y_j$ the representation for $y = g(x)$ on $[x_i, x_j]$ obtained by restricting $y_1y_2 \cdots y_k$ to $[x_i, x_j]$.  For convenience, we will also call every $y_i$ in $y_1y_2 \cdots y_k$ a node.  If $y_i = y_{i+1}$ for some $i$ (that is, $f$ is constant on $[x_i, x_{i+1}]$), we will simply write $y_1 \cdots y_iy_{i+2} \cdots y_k$ instead of 
$y_1 \cdots y_iy_{i+1}y_{i+2} \cdots y_k$.  Therefore, every two consecutive nodes in a (symbolic) representation are distinct.  

Now assume that $\{ (x_i, y_i) : i = 1, 2, \cdots, k \}$ is a set of nodes for $y = g(x)$ and $a_1a_2 \cdots a_r$ is a representation for $y = g(x)$ with $\{ a_1, a_2, \cdots, a_r \} \subset \{ y_1, y_2, \cdots, y_k \}$ and $a_i \ne a_{i+1}$ for all $1 \le i \le r-1$.  If $\{ y_1, y_2, \cdots, y_k \} \subset \{ x_1, x_2, \cdots, x_k \}$, then there is an easy way to obtain a representation for $y = g^2(x)$ from the one $a_1a_2 \cdots a_r$ for $y = g(x)$.  The procedure is as follows.  First, for any two distinct real numbers $u$ and $v$, let $[u : v]$ denote the closed interval with endpoints $u$ and $v$.  Then let $b_{i,1}b_{i,2} \cdots b_{i,t_i}$ be the representation for $y = g(x)$ on $[a_i:a_{i+1}]$ which is obtained by restricting $a_1a_2 \cdots a_r$ to $[a_i:a_{i+1}]$.  We use the following notation to indicate this fact: $a_ia_{i+1} \ra b_{i,1}b_{i,2} \cdots b_{i,t_i}$ (under $g$) if $a_i < a_{i+1}$, or $a_ia_{i+1} \ra b_{i,t_i}, \cdots, b_{i,2}b_{i,1}$ (under $g$) if $a_i > a_{i+1}$.  The above representation on $[a_i:a_{i+1}]$ exists since $\{ a_1, a_2, \cdots, a_r \} \subset \{ x_1, x_2, \cdots, x_k \}$.  Finally, if $a_i < a_{i+1}$, let $z_{i,j} = b_{i,j}$ for all $1 \le j \le t_i$.  If $a_i > a_{i+1}$, let $z_{i,j} = b_{i,t_i+1-j}$ for all $1 \le j \le t_i$.  It is easy to see that $z_{i,t_i} = z_{i+1,1}$ for all $1 \le i \le r-1$.  So, if we define $$Z = z_{1,1} \cdots z_{1,t_1}z_{2,2} \cdots z_{2,t_2} \cdots z_{r,2} \cdots z_{r,t_r},$$ then it is obvious that $Z$ is a representation for $y = g^2(x)$.  It is also obvious that the above procedure can be applied to the representation $Z$ for $y = g^2(x)$ to obtain one for $y = g^3(x)$, and so on.  

\section{Proof of Theorem 1}

In this section we let $f(x)$ denote the map as defined in Theorem 1, that is, $f(x) = -2x + 5$ if $1 \le x \le 2$, and $f(x) = x - 1$ if $2 \le x \le 3$.  The proof of part {\it (a)} of Theorem 1 will follow from two easy lemmas.  

\noindent
{\bf Lemma 4.}
{\it Under $f$, we have
$$
\begin{cases}
13 \ra 312 &, \quad 31 \ra 213, \\
12 \ra 31 &, \quad 21 \ra 13. \\
\end{cases}
$$}

In the following when we say the representation for $y = f^k(x)$, we mean the representation obtained, following the procedure as described in Section 3, by applying Lemma 4 to the representation $312$ for $y = f(x)$ successively until we get to the one for $y = f^k(x)$.

For every positive integer $k$, let $u_{1,k}$ ($u_{2,k}$ respectively) denote the number of $13$'s and $31$'s in the representation for $y = f^k(x)$ whose corresponding $x$-coordinates are $\le \, (\ge$ \, respectively) 2.  We also let $v_{1,k}$ ($v_{2,k}$ respectively) denote the number of $12$'s and $21$'s in the representation for $y = f^k(x)$ whose corresponding $x$-coordinates are $\le \, (\ge$ \, respectively) 2.  It is clear that $u_{1,1} = 1 = v_{2,1}$ and $u_{2,1} = 0 = v_{1,1}$.  Now from Lemma 4, we have

\noindent
{\bf Lemma 5.}
{\it For every positive integer $k$ and integers $i = 1, 2$, $u_{i,k+1} = u_{i,k} + v_{i,k}$ and $v_{i,k+1} = u_{i,k}$.  Furthermore, if $w_k = u_{1,k} + v_{1,k} + u_{2,k}$, then $w_1 = 1$, $w_2 = 3$, and $w_{k+2} = w_{k+1} + w_k$.  That is, $<w_k>$ is the Lucas sequence.}  

Since, for every positive integer $k$, the number of distinct solutions of the equation $f^k(x) = x$ equals $w_k$, part {\it (a)} of Theorem 1 follows from Lemma 5.  Part {\it (b)} follows from the standard inclusion-exclusion argument.  As for part {\it (c)}, we note that, for every positive integer $k$, $$a_{k+2} = 3 + \sum_{i=1}^k a_i.$$  So, for $k \ge 6$, $$(k+2)\Phi_1(k+3) > (k+2)(a_{k+3} - a_{[(k+3)/2]+1}) > (k+3)(a_{k+2} + a_{[(k+3)/2]+1}) > (k+3)\Phi_1(k+2),$$ where $[(k+3)/2]$ is the largest integer less than or equal to $(k+3)/2$.  The proof of the other statement of part {\it (c)} is easy and omitted.  This completes the proof of Theorem 1.

\section{Proof of Theorem 2}
In this section we fix any integer $n > 1$ and let $f_n(x)$ denote the map as defined in Theorem 2.  For convenience, we also let $S_n$ denote the set of all these $4n$ symbolic pairs: $i(i+1), (i+1)i$, $1 \le i \le n-1$; $n(n+2), (n+2)n, (n+1)(2n+1), (2n+1)(n+1), j(j+1), (j+1)j$, $n+2 \le j \le 2n$.  

The following lemma is easy.  

\noindent
{\bf Lemma 6.}
{\it Under $f_n$, we have
$$
\begin{cases}
n(n+2) & \ra (n+3)(n+2)n, \quad (n+2)n \ra n(n+2)(n+3), \\
(n+1)(2n+1) & \ra (n+2)n(n-1)(n-2) \cdots 321, \\
(2n+1)(n+1) & \ra 123 \cdots (n-2)(n-1)n(n+2), \\
\end{cases}
$$ and $uv \ra f_n(u)f_n(v)$ for every $uv$ in $$S_n - \{ n(n+2), (n+2)n, (n+1)(2n+1), (2n+1)(n+1) \}.$$}

In the following when we say the representation for $y = f_n^k(x)$, we mean the representation obtained, following the procedure as described in Section 3, by applying Lemma 6 to the representation $$(n+1)(2n+1)(2n)(2n-1) \cdots (n+2)n(n-1)(n-2) \cdots 321$$ for $y = f_n(x)$ successively until we get to the one for $y = f_n^k(x)$.

For every positive integer $k$ and all integers $i, j$ with $1 \le i, j \le 2n$, let $b_{k,i,j,n}$ denote the number of $uv$'s and $vu$'s in the representation for $y = f_n^k(x)$ whose corresponding $x$-coordinates are in $[i, i+1]$, where $uv = j(j+1)$ if $1 \le j \le n-1$ or $n+2 \le j \le 2n$, $uv = n(n+2)$ if $j = n$, and $uv = (n+1)(2n+1)$ if $j = n+ 1$.  It is obvious that $b_{1,1,n+1,n} = 1$, $b_{1,i,2n+2-i,n} = 1$ if $2 \le i \le n$, $b_{1,i,2n+1-i,n} = 1$ if $n+1 \le i \le 2n$, and $b_{1,i,j,n} = 0$ elsewhere.  From Lemma 6, we see that the sequences $<b_{k,i,j,n}>$ are exactly the same as those defined in Section 2.  

Since $$c_{k,n} = \sum_{i=1}^{2n} b_{k,i,i,n} + b_{k,n+1,n,n} + \sum_{i=n+2}^{2n} b_{k,i,n+1,n},$$ it is clear that $c_{k,n}$ is the number of intersection points of the graph of $y = f_n^k(x)$ with the diagonal $y = x$.  This proves part {\it (a)} of Theorem 2.  Part {\it (b)} follows from the standard inclusion-exclusion argument.  As for part {\it (c)}, we note that there exist $2n+1$ nonzero constants $\alpha_j$'s such that $$b_{k,1,n,n} = \sum_{j=1}^{2n+1} \alpha_jx_j^k$$ for all positive integers $k$, where $\{ x_j : 1 \le j \le 2n+1 \}$ is the set of all zeros (including complex zeros) of the polynomial $x^{2n+1} - 2x^{2n-1} - 1$.  Since $c_{k+2n-2,n}$ can also be expressed as $$b_{k+2n-2,n+1,n,n} + 2nb_{k,1,n,n} + \sum_{i=2}^n (2^i-2)b_{k,1,n+1-i,n},$$ Part {\it (c)} follows from property (i) of the sequences $< b_{k,i,j,n} >$ stated in Section 2.  This completes the proof of Theorem 2.

\end{document}